\documentclass[12pt]{amsart}
\usepackage{amsmath,amssymb,latexsym,cancel,rotating}
\usepackage{graphicx,amssymb,mathrsfs,amsmath,color,fancyhdr,amsthm}
\usepackage[all]{xy}
\usepackage{pgfplots}
\usepackage{tikz}
\usepackage{cite}
\usepackage{circuitikz-1.0}
\usepackage{verbatim}
\usepackage{float}
\usetikzlibrary{decorations.pathreplacing,decorations.markings}
\usetikzlibrary{arrows, decorations.pathmorphing,backgrounds,positioning,fit,petri}
\usetikzlibrary{decorations.pathreplacing,decorations.markings}
\usetikzlibrary{arrows.meta}

\newtheorem{theorem}{Theorem}[section]
\newtheorem{lemma}[theorem]{Lemma}

\newtheorem{definition}[theorem]{Definition}
\newtheorem{proposition}[theorem]{Proposition}

\newtheorem{remark}[theorem]{Remark}

\begin{document}

\title[Reflection Functors For Continuous Quivers Of Type $A$]
{Reflection Functors For Continuous Quivers Of Type $A$}

\author{Yanxiu Liu}             
\address{School of Science, Beijing Forestry University, Beijing 100083, P. R. China}
\email{liuyanxiu@bjfu.edu.cn (Y.Liu)} 

\author{Minghui ZHAO}     
\address{School of Science, Beijing Forestry University, Beijing 100083, P. R. China}
\email{zhaomh@bjfu.edu.cn (M.Zhao)}



\date{\today}

 
\bibliographystyle{abbrv}

\begin{abstract}
  As generalizations of quivers of type $A$, Igusa-Rock-Todorov in  \cite{8} introduced continuous quivers of type $A$. In this paper, we shall generalize BGP reflection functors to continuous quivers of type $A$.
  
\end{abstract}

\maketitle

\section{Introduction}

In \cite{5}, Gabriel gave the classification of indecomposable representations of a finite type quiver.
In \cite{7}, Bernstein, Gelfand and Ponomarev introduced reflection functors and gave a new proof of Gabriel's theorem.

Representations of quivers of type $A$ play an important role in persistent homology, which have been widely used in topological data analysis in \cite{13}.
In \cite{12}, Carlsson and de Silva introduced zigzag persistent homology. They also introduced diamond principle as a calculational tool, which has a direct connection with BGP reflection functors in \cite{13}.

As a generalization of Gabriel's theorem for quivers of type $A$, Crawley-Boevey in \cite{14}  and Botnan in \cite{15}  gave the classification of indecomposable representations of $\mathbb{R}$ and infinite zigzag, respectively.

In \cite{8}, Igusa-Rock-Todorov introduced continuous quivers of type $A$, which is a generalization of quivers of type $A$, $\mathbb{R}$ and infinite zigzag. They also classified indecomposable representations and proved a decomposable theorem.

In this paper, we shall introduce reflection functors for the continuous quivers of type $A$, as generalizations of BGP reflection functors and diamond principle.
  
In Section 2, we shall recall basic notations for BGP reflection functors.
Reflection functors for continuous quivers of type $A$ will be introduced in Section 3. In Section 4, we shall study the properties of this functor and give the main result.

\section{BGP reflection functors}

\subsection{Quivers}
 A quiver $Q=(Q_0,Q_1,s,t)$  is a quadruple, where: 
\begin{enumerate}
  \item[(a)]$Q_0$ is the set of vertices and $Q_1$ is the set of arrows;
  \item[(b)] $s, t: Q_1 \to Q_0$ are two maps such that $s(\alpha)$ is the source and $t(\alpha)$ is the target of $\alpha$ for any $\alpha$  in $Q_1$.
\end{enumerate}

\subsection{Representations of quivers}
Let $Q$ be a finite quiver and fix an algebraically closed field $K$.  A representation of $Q$ is $M=(M_a,\varphi _\alpha )_{a\in Q_0,  \alpha \in Q_1}$, where:
\begin{enumerate}
  \item[(a)]  $M_a$ is a $K$-vector space  for any vertex $a$ in $Q_0$;
  \item[(b)] $\varphi _\alpha$: $M_{s(\alpha)}\rightarrow M_{t(\alpha)}$ is a $K$-linear map for any arrow $\alpha: s(\alpha )\rightarrow t(\alpha )$ in $Q_1$.
\end{enumerate}

Let $M=(M_a,\varphi _\alpha )$ and $M^{'}=(M^{'}_{a},\varphi^{'} _\alpha )$ be two representations of $Q$. A $K$-linear map $f=(f_i)_{i\in{Q_0}}: M\rightarrow M^{'}$ is called a homomorphism of representations, if the following diagram is commutative
$$\xymatrix{
 M_{s(\alpha)}\ar[d]^-{f_{s(\alpha)}}\ar[r]^-{\varphi _\alpha}  &M_{t(\alpha)} \ar[d]^-{f_{t(\alpha)}}\\
  M^{'}_{s(\alpha)}\ar[r]^-{\varphi ^{'}_\alpha}&{M^{'}_{t(\alpha)}.}}$$

And the category of finite dimensional representations of $Q$ over $K$ is denoted by $rep_{K}(Q)$.

\subsection{BGP reflection functors}
Let $Q=(Q_0,Q_1,s,t)$ be a finite quiver.
The vertex $i\in Q_0$ is called a sink, if $s(\alpha)\neq i$ for any $\alpha \in Q_1$.

Consider a new quiver $\sigma _{i}Q=(Q_0,Q^{'}_1,s,t)$ by reversing the direction of arrows $\alpha$ such that $s(\alpha)=i$ or $t(\alpha)=i$.

  For instance, if $Q$ is the quiver (1),  then $\sigma _{2}Q$ is the quiver (2).
   \begin{equation}
   \quad \bullet \longrightarrow \bullet \longleftarrow \bullet \longleftarrow \bullet \longleftarrow \bullet 
    \end{equation}
  $$\qquad 1 \qquad 2\qquad 3 \qquad 4 \qquad 5 \quad$$

     \begin{equation}
   \quad \bullet \longleftarrow \bullet \longrightarrow \bullet \longleftarrow \bullet \longleftarrow \bullet
   \end{equation}
  $$\qquad 1 \qquad 2\qquad 3 \qquad 4 \qquad 5 \quad$$
  
 Next, we shall recall the definition of reflection functor $\mathsf{S}^{+}_{i} : rep_{K}(Q)\longrightarrow rep_{K}(Q^{'})$ for any $i \in Q_0 $.
 
 For any $M=(M_a,\varphi _\alpha ) \in rep_{K}(Q)$, define $\mathsf{S}^{+}_{i}M=M^{'}=(M^{'}_{a},\varphi^{'} _\alpha )$ as follows. Let 
 $$M^{'}_{a}=M_{a}\,\,\,\,\,\textrm{for $i\neq p$},$$
 and $M^{'}_{i}$ be the kernel of $$\bigoplus_ {\alpha, t(\alpha)=i}M_{s(\alpha)}\longrightarrow M_{i}, $$
 that is, we have the following exact sequence of vector spaces
 $$\xymatrix{0\ar[r]&M^{'}_{i}\ar[r]^-{\varphi _\alpha^{'}}&\bigoplus_ {\alpha, t(\alpha)=i}M_{s(\alpha)} \ar[r]^-{\varphi _\alpha }&M_{i}}.$$

Let $f=(f_{a})_{a \in Q_0}:  M\rightarrow N$ be a morphism in $rep_{K}(Q)$, where $M=(M_a,\varphi _\alpha )$ and $N=(N_a,\psi _\alpha )$. There exists a morphism
$$\mathsf{S} ^{+}_{i}(f) = f^{'} = (f^{'}_{a})_{a \in Q^{'}_0} : \mathsf{S} ^{+}_{i}M \to \mathsf{S} ^{+}_{i}N  $$
 in $rep_{K}(Q^{'})$ defined as follows.
For any $a\neq i$, let $$f^{'}_{a}=f_{a},$$  whereas $f^{'}_{i}$ is the unique morphism, making the following diagram commutative
    $$\xymatrix{
      0\ar[r]&M^{'}_{i}\ar@{-->}[d]^-{\exists f^{'}_{i}}\ar[r]&\bigoplus_{\alpha }M_{s(\alpha)} \ar[d]^-{\bigoplus_{f_{s(\alpha)}}}\ar[r]&M_i\ar[d]^-{f_i}\\
    0\ar[r]&N^{'}_{i}\ar[r]&\bigoplus_{\alpha }N_{s(\alpha)}\ar[r]&N_i.
    }$$
    
And Bernstein, Gelfand and Ponomarev studied the properties of reflection functors $\mathsf{S}^{+}_{i}$ in \cite{7}.

\section{Reflection Functors For Continuous Quivers Of Type $A$ }

In this section, firstly, we shall recall the definition of continuous quivers in \cite{8}.

\subsection{Continuous quivers of type $A$}
  A quiver of continuous type $A$, denoted by $A_\mathbb{R} $, is a triple $(\mathbb{R}, \mathbf{S}, \preccurlyeq )$, where:
  \begin{enumerate}
  \item[(a)] $\mathbf{S}\subset \mathbb{R}$ is a discrete subset;
  \item[(b)] elements of $\mathbf{S} \cup {\pm \propto}$ are indexed by a subset of $\mathbb{Z} \cup {\pm \propto}$;
  \item[(c)] the partial order $\preccurlyeq$ on $\mathbb{R}$ (which we call the orientation of $A_\mathbb{R} $)  does not change between consecutive elements of $\mathbf{S} \cup {\pm \propto}$.
  \end{enumerate}
  
  The element $S_{i}$ is called a sink, if $S_{i}\prec S_{i+1}$ and $S_{i}\prec S_{i-1}$, and ${S}_{i}$ is called a source, if $S_{i+1}\prec S_{i}$ and $S_{i-1}\prec S_{i}$ for any  $i \in \mathbb{Z} $. 
  
 \begin{remark}
There is a minute difference for the definitions of continuous quivers between this paper and  \cite{8}. In  \cite{8}, $\mathbf{S}$ is regularly punctuated with sinks and sources, that is, sinks are surrounded only by sources. In this paper, we don't specify the places of sinks and sources.
\end{remark}

\subsection{Representations of continuous quivers}
Let $A_\mathbb{R}=(\mathbb{R}, \mathbf{S}, \preccurlyeq ) $ be a quiver of continuous type $A$. And $V= (V(x), V(x, y))$ is called a representation of $A_\mathbb{R}$, if $V(x)$ is a vector space for any $x \in \mathbb{R}$, and $V(x,y):V(x) \rightarrow V(y)$ is a linear map for any $x \preccurlyeq y$ satisfied $V(y,z) \circ V(x,y)=V(x,z)$  for any $x \preccurlyeq y \preccurlyeq z$.

 Consider linear map $f_{x}:V(x)\rightarrow W(x)$ for any $x\in \mathbb{R} $, making the following diagram commutative
   $$\xymatrix{
    V(x)\ar[d]^-{f_x}\ar[r]^-{V(x,y)}&V(y) \ar[d]^-{f_y}\\
    W(x)\ar[r]^-{W(x,y)}&{W(y)}.}$$ 
    
 The collection $f= {(f_x)}_{x \in \mathbb{R}}: V\rightarrow W$ is called a morphism of representations of $A_\mathbb{R}$. And we denote the category of finite dimensional $K$-linear representations of $A_\mathbb{R}$ by $rep_{K}(A_\mathbb{R})$.

\subsection{Reflection functors}
In this subsection, we shall give the definition of reflection functors for continuous quivers of type $A$. 

Let $A_\mathbb{R}= (\mathbb{R}, \mathbf{S},\preccurlyeq )$ and $S_{k-1}$, $S_{k}$, $S_{k+1} \in {\mathbf{S}}$ be three points  satisfied $S_{k}\prec S_{k-1}$, $S_{k}\prec S_{k+1}$, and consider a new continuous quiver $\sigma _{k}A_{\mathbb{R}} = A_{\mathbb{R}^{'} } = (\mathbb{R}, {\mathbf{S}^{'}},\preccurlyeq^{'} )$ defined as follows.
\begin{enumerate}
  \item[(a)] The set $\mathbf{S}^{'}\backslash\{ {S^{'}_{k}}\}  = \mathbf{S} \backslash \{ {S_k}$\}  and $S^{'}_{k}$ satisfies that $S^{'}_{k} + S_{k} = S_{k+1} + S_{k-1}$.
  \item[(b)] The partial order $\preccurlyeq ^{'}$ on $\mathbb{R} $, which we call the orientation of $A_\mathbb{R}$, is defined as follows.
 
   For any $x \leq y \in (S^{'}_{k-1}=S_{k-1}, S^{'}_{k})$, define $x \preccurlyeq^{'} y$, and for any $x \leq y \in (S^{'}_{k},S^{'}_{k+1}=S_{k+1}, )$, define $y \preccurlyeq^{'} x$.
  \end{enumerate}

  For instance, if $A_\mathbb{R} $ is the quiver (3), then $A^{'}_\mathbb{R} =\sigma _{k}A_\mathbb{R} $ is the quiver (4).

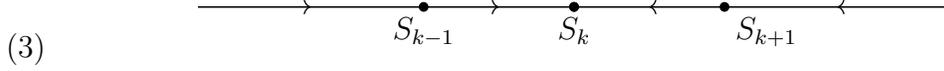
\begin{equation}
    \begin{tikzpicture}
    \begin{scope}[very thick,decoration={markings, mark=at position 0.5 with {\arrow{>}}}] 
      \draw [line width=0.5pt](5,0) to  (2,0)node[below right ]{$S_{k+1} $};
      \draw [line width=0.5pt](2,0) to [short,*-*] (0,0)node[below]{$S_{k}$};
      \draw [line width=0.5pt](0,0) to [short,*-*] (-2,0)node[below]{$S_{k-1}$};
      \draw [line width=0.5pt](-2,0) to  (-5,0);
      \draw[line width=0.5pt][postaction={decorate}] (-5,0)--(-2,0);
      \draw[line width=0.5pt][postaction={decorate}] (5,0)--(2,0);
      \draw[line width=0.5pt][postaction={decorate}] (2,0)--(0,0);
      \draw[line width=0.5pt][postaction={decorate}] (-2,0)--(0,0);
    \end{scope}
  \end{tikzpicture}
\end{equation}

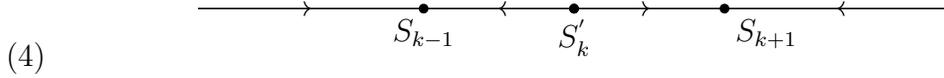
\begin{equation}
    \begin{tikzpicture}
    \begin{scope}[very thick,decoration={markings, mark=at position 0.5 with {\arrow{>}}}] 
      \draw [line width=0.5pt](5,0) to (2,0)node[below right]{$S_{k+1}$};
      \draw [line width=0.5pt](2,0) to [short,*-*] (0,0)node[below]{$S^{'}_{k}$};
      \draw [line width=0.5pt](0,0) to [short,*-*] (-2,0)node[below]{$S_{k-1}$};
      \draw [line width=0.5pt](-2,0) to (-5,0);
      \draw[line width=0.5pt][postaction={decorate}] (-5,0)--(-2,0);
      \draw[line width=0.5pt][postaction={decorate}] (5,0)--(2,0);
      \draw[line width=0.5pt][postaction={decorate}] (0,0)--(2,0);
      \draw[line width=0.5pt][postaction={decorate}] (0,0)--(-2,0);
    \end{scope}
  \end{tikzpicture}
\end{equation}

If we put $A_\mathbb{R}$ and $A^{'}_\mathbb{R}$ together, we get the following diagram.

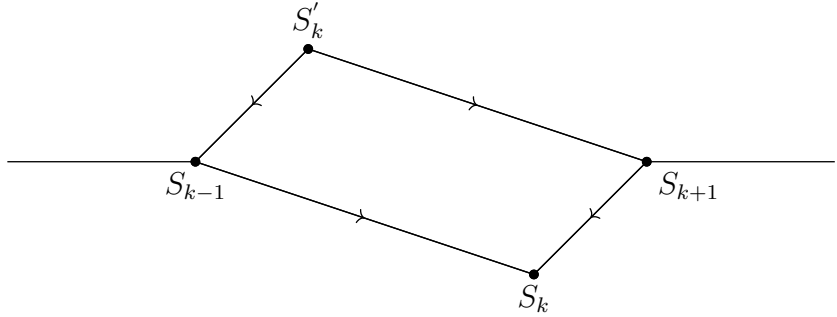
\begin{center}
  \begin{tikzpicture}
  \begin{scope}[very thick,decoration={markings, mark=at position 0.5 with {\arrow{>}}}] 
    \draw [line width=0.5pt](5.5,0) to (3,0)node[below right]{$S_{k+1}$};
    \draw [line width=0.5pt](3,0) to [short,*-*](-1.5,1.5)node[above]{$S^{'}_{k}$};
    \draw [line width=0.5pt](3,0) to  [short,*-*] (1.5,-1.5)node[below]{$S_{k}$};
    \draw [line width=0.5pt](-1.5,1.5) to [short,*-*] (-3,0)node[below]{$S_{k-1}$};
    \draw [line width=0.5pt](1.5,-1.5) to [short,*-*] (-3,0);
    \draw [line width=0.5pt](-5.5,0) to (-3,0);
    \draw[line width=0.5pt][postaction={decorate}] (-1.5,1.5)--(-3,0);
    \draw[line width=0.5pt][postaction={decorate}] (-1.5,1.5)--(3,0);
    \draw[line width=0.5pt][postaction={decorate}] (-3,0)--(1.5,-1.5);
    \draw[line width=0.5pt][postaction={decorate}] (3,0)--(1.5,-1.5);
  \end{scope}
\end{tikzpicture}
\end{center}

  Now, we shall give the definition of reflection functors. Let  $S_k$ be a sink, then we shall define the reflection functor
  $$\mathsf{S} _{k}^{+}: rep_{K}(A_\mathbb{R})\longrightarrow rep_{K}(A^{'}_\mathbb{R}).$$

For any object $V \in rep_{K}(A_\mathbb{R})$, we can give the  definition of representation $\mathsf{S}^{+}_{k}V=V^{'}$. If $x\notin ({S_{k-1},{S_{k+1}}})$, let $V^{'}(x)=V(x)$. If  $x\in ({S_{k-1},{S_{k+1}}})$, let $x^{'}=S_{k+1}+S_{k-1}-x$, and define $V^{'}(x)$ as follows.

  \begin{enumerate}
    \item[(a)]  For any $x\in ({S_{k-1},{S^{'}_{k}}})$,

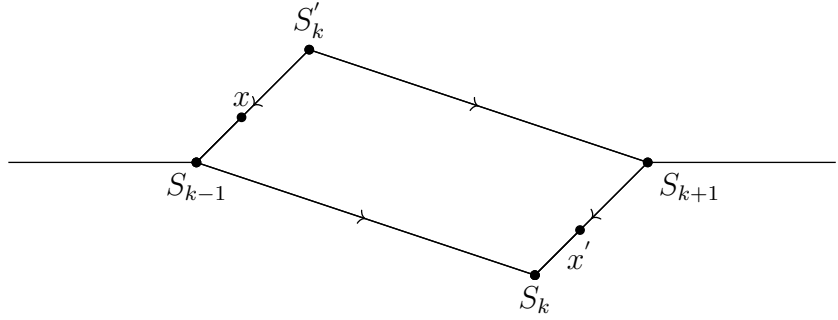
    \begin{center}
  \begin{tikzpicture}
  \begin{scope}[very thick,decoration={markings, mark=at position 0.5 with {\arrow{>}}}] 
    \draw [line width=0.5pt](5.5,0) to (3,0)node[below right]{$S_{k+1}$};
    \draw [line width=0.5pt](-3,0) to [short,*-*] (-2.4,0.6)node[above]{$x$};
    \draw [line width=0.5pt](3,0) to [short,*-*] (-1.5,1.5)node[above]{$S^{'}_{k}$};
    \draw [line width=0.5pt](3,0) to [short,*-*] (1.5,-1.5)node[below]{$S_{k}$};
    \draw [line width=0.5pt](-1.5,1.5) to [short,*-*] (-3,0)node[below]{$S_{k-1}$};
    \draw [line width=0.5pt](1.5,-1.5) to [short,*-*] (2.1,-0.9)node[below]{$x^{'}$};
    \draw [line width=0.5pt](1.5,-1.5) to [short,*-*] (-3,0);
    \draw [line width=0.5pt](-5.5,0) to (-3,0);
    \draw[line width=0.5pt][postaction={decorate}] (-1.5,1.5)--(-3,0);
    \draw[line width=0.5pt][postaction={decorate}] (-1.5,1.5)--(3,0);
    \draw[line width=0.5pt][postaction={decorate}] (-3,0)--(1.5,-1.5);
    \draw[line width=0.5pt][postaction={decorate}] (3,0)--(1.5,-1.5);
  \end{scope}
\end{tikzpicture}
\end{center}
 $V^{'}(x)$ is defined as the kernel of
    $$V(S_{k-1})\oplus V(x^{'}) \longrightarrow V(S_{k}),$$
    that is, we have the following short exact sequence of vector spaces:
    $$0 \longrightarrow V^{'}(x) \longrightarrow V(S_{k-1})\oplus V(x^{'}) \overset{D_{1}}{\longrightarrow} V(S_{k}),$$
    where $D_{1}=(V(S_{k-1},S_{k}), -V(x^{'},S_{k}))$.\\
    
    Denote the map from $V^{'}(x) =\mathsf{S} _{k}^{+} V(x)$ to $V(x^{'})$ in the exact sequence by $C_{V}(x)$. In this case, the following diagram is a pull-back
    
       \xymatrix{
    &&&&& V^{'}(x)\ar[d]^{V^{'}(x, S_{k-1})} \ar[rr] ^{C_{V}(x)}&&V(x^{'}) \ar[d]^{V(x^{'}, S_{k})}\\
    &&&&& V(S_{k-1})\ar[rr] ^{V(S_{k-1},S_{k} )}&&{V(S_{k}).}
    }

    \item[(b)]  For any $x\in [{S^{'}_{k},{S_{k+1}}})$, 
  
    
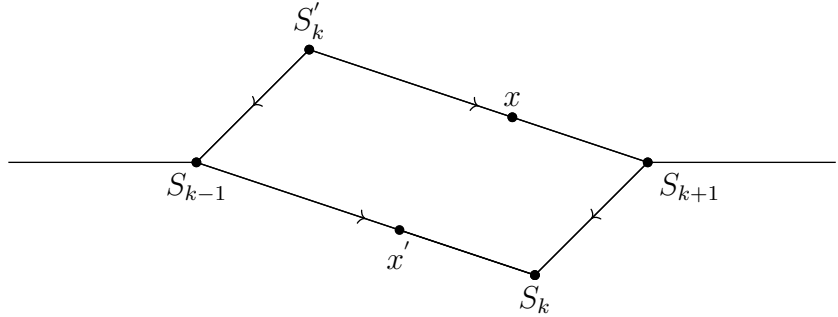
\begin{center}
  \begin{tikzpicture}
  \begin{scope}[very thick,decoration={markings, mark=at position 0.5 with {\arrow{>}}}] 
    \draw [line width=0.5pt](5.5,0) to (3,0)node[below right]{$S_{k+1}$};
    \draw [line width=0.5pt](3,0) to [short,*-*] (1.2,0.6)node[above]{$x$};
    \draw [line width=0.5pt](1.2,0.6) to [short,*-*] (-1.5,1.5)node[above]{$S^{'}_{k}$};
    \draw [line width=0.5pt](3,0) to [short,*-*] (1.5,-1.5)node[below]{$S_{k}$};
    \draw [line width=0.5pt](-1.5,1.5) to [short,*-*] (-3,0)node[below]{$S_{k-1}$};
    \draw [line width=0.5pt](1.5,-1.5) to [short,*-*] (-0.3,-0.9)node[below]{$x^{'}$};
    \draw [line width=0.5pt](1.5,-1.5) to [short,*-*] (-3,0);
    \draw [line width=0.5pt](-5.5,0) to (-3,0);
    \draw[line width=0.5pt][postaction={decorate}] (-1.5,1.5)--(-3,0);
    \draw[line width=0.5pt][postaction={decorate}] (-1.5,1.5)--(3,0);
    \draw[line width=0.5pt][postaction={decorate}] (-3,0)--(1.5,-1.5);
    \draw[line width=0.5pt][postaction={decorate}] (3,0)--(1.5,-1.5);
  \end{scope}
\end{tikzpicture}
\end{center}
 $V^{'}(x)$ is defined as the kernel of
  $$V(x^{'})\oplus V(S_{k+1}) \longrightarrow V(S_{k}),$$
  that is, we have the following short exact sequence of vector spaces:
  $$0 \longrightarrow V^{'}(x) \longrightarrow V(x^{'})\oplus V(S_{k+1}) \overset{D_{2}}{\longrightarrow} V(S_{k}),$$
  where $D_{2}=(V(x^{'},S_{k}),- V(S_{k-1},S_{k}))$.\\

    Denote the map from $V^{'}(x) =\mathsf{S} _{k}^{+} V(x)$ to $V(x^{'})$ in the exact sequence by $C_{V}(x)$. In this case, the following diagram is a pull-back    
    
  \xymatrix{
    &&&&& V^{'}(S^{'}_{k})\ar[rr]^{V^{'}(S^{'}_{k}, S_{k+1})} \ar[d] ^{V^{'}(S^{'}_{k}, x)}&& V(S_{k+1})\ar[d] ^{V(S_{k+1}, x^{'})}\\
    &&&&& V^{'}(x) \ar[rr] ^{C_{V}(x)}&&{V(x^{'}).}\\
    }

    \end{enumerate}

For any $V \in rep_{K}(A_\mathbb{R})$, if $x \preccurlyeq y\notin ({S_{k-1},{S_{k+1}}})$, let $V^{'}(x,y)=V(x,y)$. If  $x\in ({S_{k-1},{S_{k+1}}})$, let $x^{'}=S_{k+1}+S_{k-1}-x$, $y^{'}=S_{k+1}+S_{k-1}-y$, and define $V^{'}(x,y)$ as follows.

   \begin{enumerate}
    \item[(a)]  Let $x,y\in ({S_{k-1},{S^{'}_{k}}})$ such that $x\prec^{'} y$. Note that $x^{'}\prec y^{'}$.

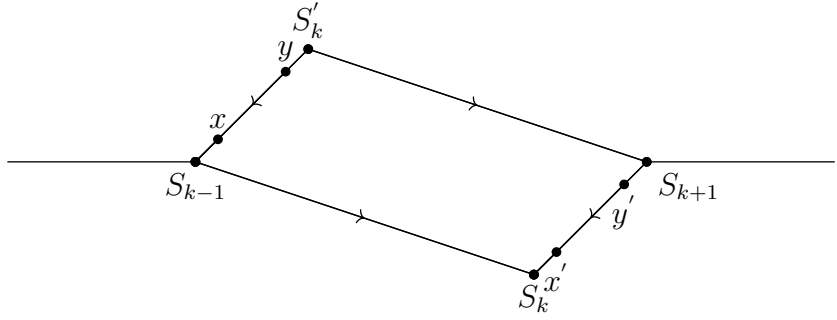
    \begin{center}
  \begin{tikzpicture}
  \begin{scope}[very thick,decoration={markings, mark=at position 0.5 with {\arrow{>}}}] 
    \draw [line width=0.5pt](5.5,0) to (3,0)node[below right]{$S_{k+1}$};
    \draw [line width=0.5pt](-3,0) to [short,*-*] (-2.7,0.3)node[above]{$x$};
     \draw [line width=0.5pt](-3,0) to [short,*-*] (-1.8,1.2)node[above]{$y$};
    \draw [line width=0.5pt](3,0) to [short,*-*] (-1.5,1.5)node[above]{$S^{'}_{k}$};
    \draw [line width=0.5pt](3,0) to [short,*-*] (1.5,-1.5)node[below]{$S_{k}$};
    \draw [line width=0.5pt](-1.5,1.5) to [short,*-*] (-3,0)node[below]{$S_{k-1}$};
    \draw [line width=0.5pt](1.5,-1.5) to [short,*-*] (1.8,-1.2)node[below]{$x^{'}$};
     \draw [line width=0.5pt](1.5,-1.5) to [short,*-*] (2.7,-0.3)node[below]{$y^{'}$};
    \draw [line width=0.5pt](1.5,-1.5) to [short,*-*] (-3,0);
    \draw [line width=0.5pt](-5.5,0) to (-3,0);
    \draw[line width=0.5pt][postaction={decorate}] (-1.5,1.5)--(-3,0);
    \draw[line width=0.5pt][postaction={decorate}] (-1.5,1.5)--(3,0);
    \draw[line width=0.5pt][postaction={decorate}] (-3,0)--(1.5,-1.5);
    \draw[line width=0.5pt][postaction={decorate}] (3,0)--(1.5,-1.5);
  \end{scope}
\end{tikzpicture}
\end{center}

 Since $V^{'}(x)$ is the kernel of $V(S_{k-1})\oplus V(x^{'}) \to V(S_{k})$, and $V^{'}(y)$ is the kernel of $V(S_{k-1})\oplus V(y^{'}) \to V(S_{k})$, there must exist a unique map $V^{'}(x, y): V^{'}(x) \to V^{'}(y)$ making the following diagram commutative
    $$\xymatrix{
    0\ar[r]&V^{'}(x)\ar[d]^{ V^{'}(x, y)}\ar[r]&V(S_{k-1}) \oplus V(x^{'}) \ar[d]^{\bigl( \begin{smallmatrix} id & \\  & {V(x^{'}, y^{'})} \end{smallmatrix} \bigr)} \ar[r]&V(S_{k})\ar[d]^{id}\\
    0\ar[r]&V^{'}(y)\ar[r]&V(S_{k-1})\oplus V(y^{'})\ar[r]&{V(S_{k}).}
     }$$

    \item[(b)]  Let $x,y\in [{S^{'}_{k},{S_{k+1}}})$ such that $x\prec^{'} y$. Note that $x^{'}\prec y^{'}$.
  
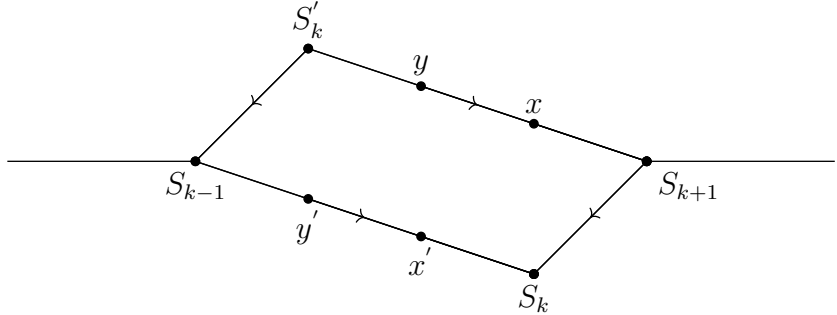
   \begin{center}
  \begin{tikzpicture}
  \begin{scope}[very thick,decoration={markings, mark=at position 0.5 with {\arrow{>}}}] 
    \draw [line width=0.5pt](5.5,0) to (3,0)node[below right]{$S_{k+1}$};
    \draw [line width=0.5pt](3,0) to [short,*-*] (1.5,0.5)node[above]{$x$};
    \draw [line width=0.5pt](3,0) to [short,*-*] (0,1)node[above]{$y$};
    \draw [line width=0.5pt](3,0) to [short,*-*] (-1.5,1.5)node[above]{$S^{'}_{k}$};
    \draw [line width=0.5pt](3,0) to [short,*-*] (1.5,-1.5)node[below]{$S_{k}$};
    \draw [line width=0.5pt](-1.5,1.5) to [short,*-*] (-3,0)node[below]{$S_{k-1}$};
    \draw [line width=0.5pt](1.5,-1.5) to [short,*-*] (0,-1)node[below]{$x^{'}$};
    \draw [line width=0.5pt](1.5,-1.5) to [short,*-*] (-1.5,-0.5)node[below]{$y^{'}$};
    \draw [line width=0.5pt](1.5,-1.5) to [short,*-*] (-3,0);
    \draw [line width=0.5pt](-5.5,0) to (-3,0);
    \draw[line width=0.5pt][postaction={decorate}] (-1.5,1.5)--(-3,0);
    \draw[line width=0.5pt][postaction={decorate}] (-1.5,1.5)--(3,0);
    \draw[line width=0.5pt][postaction={decorate}] (-3,0)--(1.5,-1.5);
    \draw[line width=0.5pt][postaction={decorate}] (3,0)--(1.5,-1.5);
  \end{scope}
\end{tikzpicture}
\end{center}

  Since $V^{'}(x)$ is the kernel of $V(x^{'})\oplus V(S_{k+1}) \to V(S_{k})$, and $V^{'}(y)$ is the kernel of $V(y^{'})\oplus V(S_{k+1}) \to V(S_{k})$, there must exist a unique map $V^{'}(x, y): V^{'}(x) \to V^{'}(y)$ making the following diagram commutative
    $$\xymatrix{
    0\ar[r]&V^{'}(x)\ar[d]^{ V^{'}(x, y)}\ar[r]&V(x^{'}) \oplus V(S_{k+1}) \ar[d]^{\bigl( \begin{smallmatrix} {V(x^{'}, y^{'})} & \\  & id \end{smallmatrix} \bigr)} \ar[r]&V(S_{k})\ar[d]^{id}\\
    0\ar[r]&V^{'}(y)\ar[r]&V(y^{'}) \oplus V(S_{k+1})  \ar[r]&{V(S_{k}).}
    }$$

    \end{enumerate}

For any morphism $f:V \to W$ in $rep_{K}(A_\mathbb{R})$, we can give the definition of a morphism $\mathsf{S}^{+}_{k}f=f^{'}=(f^{'}_{x}): \mathsf{S}^{+}_{k}V\longrightarrow \mathsf{S}^{+}_{k}W$. If $x\notin ({S_{k-1},{S_{k+1}}})$, let $f^{'}_{x}=f_{x}$.  If  $x\in ({S_{k-1},{S_{k+1}}})$, we can define  $f^{'}_{x}$ as follows. Denote $V^{'}=\mathsf{S}_{k}^{+}V$, $W^{'}=\mathsf{S}_{k}^{+}W$.
    \begin{enumerate}
      \item[(a)] Fix $x\in ({S_{k-1},{S^{'}_{k}}})$, since $V^{'}(x)$ is the kernel of $V(S_{k-1})\oplus V(x^{'}) \to V(S_{k})$, and $W^{'}(x)$ is the kernel of $W(S_{k-1})\oplus W(x^{'}) \to W(S_{k})$, there must be a map $f^{'}_{x}$ making the following diagram commutative
                  $$\xymatrix{
                  0\ar[r]&V^{'}(x)\ar[d]^-{ f^{'}_{x}}\ar[r]&V(S_{k-1})   \oplus V(x^{'}) \ar[d]^{\bigl( \begin{smallmatrix} {f_{S_{k-1}}} & \\  & {{f_{x^{'}}}} \end{smallmatrix} \bigr)} \ar[r]&V(S_{k})\ar[d]^{f_{S_k}}  \\
                  0\ar[r]&W^{'}(x)\ar[r]&W(S_{k-1})\oplus W(x^{'})\ar[r]&{W(S_{k}).}
                  }$$
         
      \item[(b)]  Similarly,  fix $x\in [{S^{'}_{k},{S_{k+1}}})$, since $V^{'}(x)$ is the kernel of $V(x^{'})\oplus V(S_{k+1}) \to V(S_{k})$, and $W^{'}(x)$ is the kernel of $W(x^{'})\oplus W(S_{k+1}) \to W(S_{k})$, there must be a map $f^{'}_{x}$ making the following diagram commutative
                  $$\xymatrix{
                  0\ar[r]&V^{'}(x)\ar[d]^-{ f^{'}_{x}}\ar[r]&V(x^{'})   \oplus V(S_{k+1}) \ar[d]^{\bigl( \begin{smallmatrix} {{f_{x^{'}}}}  & \\  &{f_{S_{k+1}}} \end{smallmatrix} \bigr)}\ar[r]&V(S_{k})\ar[d]^{f_{S_k}}\\
                  0\ar[r]&W^{'}(x)\ar[r]&W(x^{'})\oplus W(S_{k+1})\ar[r]&{W(S_{k}).}
                  }$$
                \end{enumerate}

    In a similar way, if $S_k^{'}$ is a source, we can give the definition of reflection functor $\mathsf{S}^{-}_k$
    $$\mathsf{S}^{-}_k:  rep(A^{'}_\mathbb{R}) \longrightarrow rep(A_\mathbb{R}),$$
     $$V^{'} \mapsto V=\mathsf{S}^{-}_k V^{'},$$  
     $$f^{'} \mapsto f=\mathsf{S}^{-}_k f^{'}.$$
     
For any object $V^{'} \in rep_{K}(A^{'} _\mathbb{R})$, we can give the definition of representation $\mathsf{S}^{-}_{k}V^{'}=V$. If $x\notin ({S_{k-1},{S_{k+1}}})$, let $V(x)=V^{'}(x)$. If  $x\in ({S_{k-1},{S_{k+1}}})$, let $x^{'}=S_{k+1}+S_{k-1}-x$, and define $V(x)$ as follows.

  \begin{enumerate}
    \item[(a)] For any $x\in ({S_{k},{S_{k+1}}})$,

    
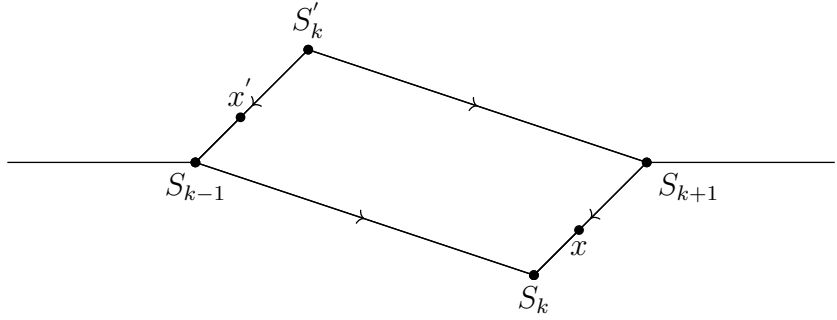
\begin{center}
  \begin{tikzpicture}
  \begin{scope}[very thick,decoration={markings, mark=at position 0.5 with {\arrow{>}}}] 
    \draw [line width=0.5pt](5.5,0) to (3,0)node[below right]{$S_{k+1}$};
    \draw [line width=0.5pt](-3,0) to [short,*-*] (-2.4,0.6)node[above]{$x^{'}$};
    \draw [line width=0.5pt](3,0) to [short,*-*] (-1.5,1.5)node[above]{$S^{'}_{k}$};
    \draw [line width=0.5pt](3,0) to [short,*-*] (1.5,-1.5)node[below]{$S_{k}$};
    \draw [line width=0.5pt](-1.5,1.5) to [short,*-*] (-3,0)node[below]{$S_{k-1}$};
    \draw [line width=0.5pt](1.5,-1.5) to [short,*-*] (2.1,-0.9)node[below]{$x$};
    \draw [line width=0.5pt](1.5,-1.5) to [short,*-*] (-3,0);
    \draw [line width=0.5pt](-5.5,0) to (-3,0);
    \draw[line width=0.5pt][postaction={decorate}] (-1.5,1.5)--(-3,0);
    \draw[line width=0.5pt][postaction={decorate}] (-1.5,1.5)--(3,0);
    \draw[line width=0.5pt][postaction={decorate}] (-3,0)--(1.5,-1.5);
    \draw[line width=0.5pt][postaction={decorate}] (3,0)--(1.5,-1.5);
  \end{scope}
\end{tikzpicture}
\end{center}
 $V(x)$ is defined as the cokernel of
    $$ V^{'}(S^{'}_{k}) \longrightarrow V^{'}(S_{k+1})\oplus V^{'}(x^{'}),$$
    that is, we have the following short exact sequence of vector spaces:
    $$ V^{'}(S^{'}_{k}) \overset{D^{'}_{1}}{\longrightarrow} V^{'}(S_{k+1})\oplus V^{'}(x^{'}) \longrightarrow  V(x)  \longrightarrow 0,$$
    where $D^{'}_{1}=\binom{V^{'}(S^{'}_{k}, S_{k+1})}{-V^{'}(S^{'}_{k}, x^{'})}$.\\
    
     Denote the map from $V^{'}(x^{'})$ to $V(x)= \mathsf{S} _{k}^{-} V^{'}(x)$ in the exact sequence by $D_{V}(x)$. In this case, the following diagram is a push-out 
     
  \xymatrix{
    &&&&& V^{'}(S^{'}_{k})\ar[rr]^{V^{'}(S^{'}_{k}, S_{k+1})} \ar[d]^{V^{'}(S^{'}_{k}, x^{'})} && V^{'}(S_{k+1})\ar[d]^{V(S_{k+1}, x)} \\
    &&&&& V^{'}(x^{'}) \ar[rr]^{D_{V}(x)} &&{V(x).}\\
    }

    \item[(b)] For any $x\in [{S_{k-1},{S_{k}}})$, 
  
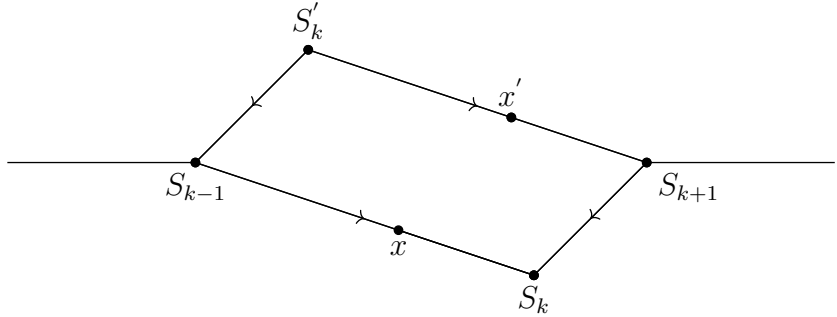
    \begin{center}
  \begin{tikzpicture}
  \begin{scope}[very thick,decoration={markings, mark=at position 0.5 with {\arrow{>}}}] 
    \draw [line width=0.5pt](5.5,0) to (3,0)node[below right]{$S_{k+1}$};
    \draw [line width=0.5pt](3,0) to [short,*-*] (1.2,0.6)node[above]{$x^{'}$};
    \draw [line width=0.5pt](1.2,0.6) to [short,*-*] (-1.5,1.5)node[above]{$S^{'}_{k}$};
    \draw [line width=0.5pt](3,0) to [short,*-*] (1.5,-1.5)node[below]{$S_{k}$};
    \draw [line width=0.5pt](-1.5,1.5) to [short,*-*] (-3,0)node[below]{$S_{k-1}$};
    \draw [line width=0.5pt](1.5,-1.5) to [short,*-*] (-0.3,-0.9)node[below]{$x$};
    \draw [line width=0.5pt](1.5,-1.5) to [short,*-*] (-3,0);
    \draw [line width=0.5pt](-5.5,0) to (-3,0);
    \draw[line width=0.5pt][postaction={decorate}] (-1.5,1.5)--(-3,0);
    \draw[line width=0.5pt][postaction={decorate}] (-1.5,1.5)--(3,0);
    \draw[line width=0.5pt][postaction={decorate}] (-3,0)--(1.5,-1.5);
    \draw[line width=0.5pt][postaction={decorate}] (3,0)--(1.5,-1.5);
  \end{scope}
\end{tikzpicture}
\end{center}
 $V(x)$ is defined as the cokernel of
    $$ V^{'}(S^{'}_{k}) \longrightarrow V^{'}(x^{'}) \oplus V^{'}(S_{k-1}),$$
    that is, we have the following short exact sequence of vector spaces:
    $$ V^{'}(S^{'}_{k}) \overset{D^{'}_{2}}{\longrightarrow}V^{'}(x^{'})  \oplus V^{'}(S_{k-1}) \longrightarrow  V(x)  \longrightarrow 0,$$
    where $D^{'}_{2}=\binom{V^{'}(S^{'}_{k}, x^{'})} {-V^{'}(S^{'}_{k}, S_{k-1})}$.\\
   
       Denote the map from $V^{'}(x^{'})$ to $V(x)= \mathsf{S} _{k}^{-} V^{'}(x)$ in the exact sequence by $D_{V}(x)$. In this case, the following diagram is a push-out 
     
              \xymatrix{
    &&&&& V^{'}(S^{'}_{k})\ar[d] ^{V^{'}(S^{'}_{k}, S_{k-1})}\ar[rr] ^{V^{'}(S^{'}_{k}, )x^{'}}&&V^{'}(x^{'}) \ar[d]^{D_{V}(x)}\\
    &&&&& V^{'}(S_{k-1})\ar[rr]^{V(S_{k-1}, x)} &&{V(x)
    }.}
 

   For any $V^{'} \in rep_{K}(A^{'}_\mathbb{R})$, if $x \preccurlyeq y\notin ({S_{k-1},{S_{k+1}}})$, let $V(x, y)=V^{'}(x, y)$. If  $x\in ({S_{k-1},{S_{k+1}}})$, let $x^{'}=S_{k+1}+S_{k-1}-x$, $y^{'}=S_{k+1}+S_{k-1}-y$, and define $V(x, y)$ as follows.

   \begin{enumerate}
    \item[(a)] Let $x,y\in ({S_{k},{S_{k+1}}})$ such that $x\prec^{'} y$. Note that $x^{'}\prec y^{'}$.

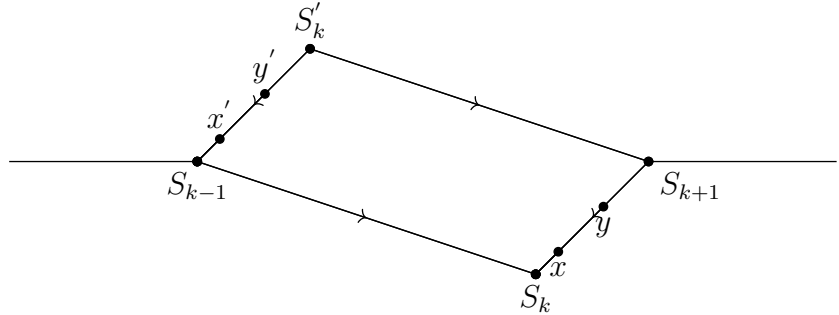
    \begin{center}
  \begin{tikzpicture}
  \begin{scope}[very thick,decoration={markings, mark=at position 0.5 with {\arrow{>}}}] 
    \draw [line width=0.5pt](5.5,0) to (3,0)node[below right]{$S_{k+1}$};
    \draw [line width=0.5pt](-3,0) to [short,*-*] (-2.7,0.3)node[above]{$x^{'}$};
     \draw [line width=0.5pt](-3,0) to [short,*-*] (-2.1,0.9)node[above]{$y^{'}$};
    \draw [line width=0.5pt](3,0) to [short,*-*] (-1.5,1.5)node[above]{$S^{'}_{k}$};
    \draw [line width=0.5pt](3,0) to [short,*-*] (1.5,-1.5)node[below]{$S_{k}$};
    \draw [line width=0.5pt](-1.5,1.5) to [short,*-*] (-3,0)node[below]{$S_{k-1}$};
    \draw [line width=0.5pt](1.5,-1.5) to [short,*-*] (1.8,-1.2)node[below]{$x$};
     \draw [line width=0.5pt](1.5,-1.5) to [short,*-*] (2.4,-0.6)node[below]{$y$};
    \draw [line width=0.5pt](1.5,-1.5) to [short,*-*] (-3,0);
    \draw [line width=0.5pt](-5.5,0) to (-3,0);
    \draw[line width=0.5pt][postaction={decorate}] (-1.5,1.5)--(-3,0);
    \draw[line width=0.5pt][postaction={decorate}] (-1.5,1.5)--(3,0);
    \draw[line width=0.5pt][postaction={decorate}] (-3,0)--(1.5,-1.5);
    \draw[line width=0.5pt][postaction={decorate}] (3,0)--(1.5,-1.5);
  \end{scope}
\end{tikzpicture}
\end{center}

          Since $V(x)$ is the cokernel of $ V^{'}(S^{'}_{k}) \longrightarrow V^{'}(S_{k+1})\oplus V^{'}(x^{'})$, and $V(y)$ is the cokernel of $ V^{'}(S^{'}_{k}) \longrightarrow V^{'}(S_{k+1})\oplus V^{'}(y^{'})$, there must exist a unique map $V(x, y): V(x) \to V(y)$ making the following diagram commutative
    $$\xymatrix{
    V^{'}(S^{'}_{k})\ar[d]^{id}\ar[r]&V^{'}(S_{k+1})\oplus V^{'}(x^{'}) \ar[d]^{\bigl( \begin{smallmatrix} id & \\  &  {V^{'}(x^{'}, y^{'})}\end{smallmatrix} \bigr)}  \ar[r]&V(x)\ar[d]^{V(x, y)}\ar[r]&0\\
    V^{'}(S^{'}_{k})\ar[r]&V^{'}(S_{k+1})\oplus V^{'}(y^{'})\ar[r]&{V(y)}\ar[r]&0.
     }$$

   
    \item[(b)]  Let $x,y\in [{S_{k-1},{S_{k}}})$ such that $x\prec^{'} y$. Note that $x^{'}\prec y^{'}$.
  
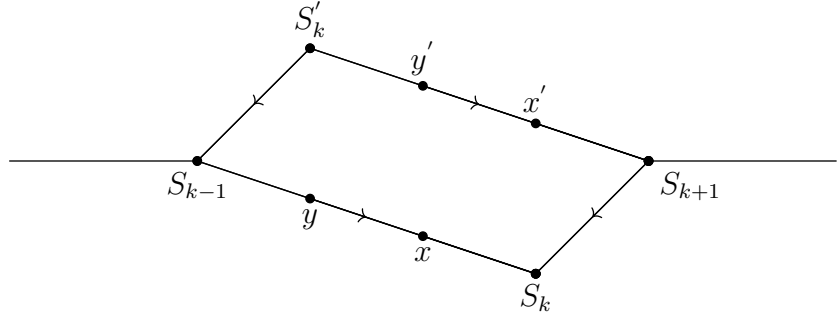
   \begin{center}
  \begin{tikzpicture}
  \begin{scope}[very thick,decoration={markings, mark=at position 0.5 with {\arrow{>}}}] 
    \draw [line width=0.5pt](5.5,0) to (3,0)node[below right]{$S_{k+1}$};
    \draw [line width=0.5pt](3,0) to [short,*-*] (1.5,0.5)node[above]{$x^{'}$};
    \draw [line width=0.5pt](3,0) to [short,*-*] (0,1)node[above]{$y^{'}$};
    \draw [line width=0.5pt](3,0) to [short,*-*] (-1.5,1.5)node[above]{$S^{'}_{k}$};
    \draw [line width=0.5pt](3,0) to [short,*-*] (1.5,-1.5)node[below]{$S_{k}$};
    \draw [line width=0.5pt](-1.5,1.5) to [short,*-*] (-3,0)node[below]{$S_{k-1}$};
    \draw [line width=0.5pt](1.5,-1.5) to [short,*-*] (0,-1)node[below]{$x$};
    \draw [line width=0.5pt](1.5,-1.5) to [short,*-*] (-1.5,-0.5)node[below]{$y$};
    \draw [line width=0.5pt](1.5,-1.5) to [short,*-*] (-3,0);
    \draw [line width=0.5pt](-5.5,0) to (-3,0);
    \draw[line width=0.5pt][postaction={decorate}] (-1.5,1.5)--(-3,0);
    \draw[line width=0.5pt][postaction={decorate}] (-1.5,1.5)--(3,0);
    \draw[line width=0.5pt][postaction={decorate}] (-3,0)--(1.5,-1.5);
    \draw[line width=0.5pt][postaction={decorate}] (3,0)--(1.5,-1.5);
  \end{scope}
\end{tikzpicture}
\end{center}

  Since $V(x)$ is the cokernel of $V^{'}(S^{'}_{k}) \longrightarrow V^{'}(x^{'}) \oplus V^{'}(S_{k-1})$, and $V(y)$ is the cokernel of $V^{'}(S^{'}_{k}) \longrightarrow V^{'}(y^{'}) \oplus V^{'}(S_{k-1})$, there must exist a unique map $V(x, y): V(x) \to V(y)$ making the following diagram commutative
    $$\xymatrix{
    V^{'}(S^{'}_{k})\ar[d]^{id}\ar[r]& V^{'}(x^{'}) \oplus V^{'}(S_{k-1}) \ar[d]^{\bigl( \begin{smallmatrix} {V^{'}(x^{'}, y^{'})} & \\  &  id\end{smallmatrix} \bigr)}  \ar[r]&V(x)\ar[d]^{V(x ,y)}\ar[r]&0\\
    V^{'}(S^{'}_{k})\ar[r]& V^{'}(y^{'}) \oplus V^{'}(S_{k-1})  \ar[r]&{V(y)}\ar[r]&0.
    }$$
      \end{enumerate}
       
For any morphism $f^{'} : V^{'} \to W^{'}$ in $rep_{K}(A^{'}_\mathbb{R})$, we can give the definition of a morphism $\mathsf{S}^{-}_{k}f^{'}=f  =(f_x): \mathsf{S}^{-}_{k}V^{'}\longrightarrow \mathsf{S}^{-}_{k}W^{'}$. If $x\notin ({S_{k-1},{S_{k+1}}})$, let $f_{x}=f^{'}_{x}$.  If  $x\in ({S_{k-1},{S_{k+1}}})$, we can define $f_{x}$ as follows. Denote $V=\mathsf{S}_{k}^{-}V^{'}$, $W=\mathsf{S}_{k}^{-}W^{'}$.
    \begin{enumerate}
      \item[(a)]  Fix $x \in (S_{k}, S_{k+1})$, since $V(x)$ is the cokernel of $ V^{'}(S^{'}_{k}) \longrightarrow V^{'}(S_{k+1})\oplus V^{'}(x^{'})$, and $W(x)$ is the cokernel of $ W^{'}(S^{'}_{k}) \longrightarrow W^{'}(S_{k+1})\oplus W^{'}(x^{'})$, there must be a map $f_{x}: V(x) \to W(x)$ making the following diagram commutative
              
     $$\xymatrix{
    V^{'}(S^{'}_{k})\ar[d]^{f^{'}_{S^{'}_{k}}}\ar[r]&V^{'}(S_{k+1}) \oplus V^{'}(x^{'}) \ar[d]^{\bigl( \begin{smallmatrix} {f^{'}_{S_{k+1}}}  & \\  &{f^{'}_{x^{'}}} \end{smallmatrix} \bigr)}\ar[r]&V(x)\ar[d]^{f_{x}}\ar[r]&0\\
    W^{'}(S^{'}_{k})\ar[r]&W^{'}(S_{k+1})\oplus W^{'}(x^{'})\ar[r]&{W(x)}\ar[r]&0.
     }$$
         
      \item[(b)]  Similarly, fix $x \in (S_{k-1}, S_{k})$, since $V(x)$ is the cokernel of $V^{'}(S^{'}_{k}) \longrightarrow V^{'}(x^{'}) \oplus V^{'}(S_{k-1})$, and $W(y)$ is the cokernel of $W^{'}(S^{'}_{k}) \longrightarrow W^{'}(y^{'}) \oplus W^{'}(S_{k-1})$, there must be a map $f_{x}: V(x) \to W(x)$ making the following diagram commutative
    $$\xymatrix{
    V^{'}(S^{'}_{k})\ar[d]^{f^{'}_{S^{'}_{k}}}\ar[r]& V^{'}(x^{'}) \oplus V^{'}(S_{k-1})\ar[d]^{\bigl( \begin{smallmatrix} {f^{'}_{x^{'}}}   & \\  &{f^{'}_{S_{k-1}}}\end{smallmatrix} \bigr)} \ar[r]&V(x)\ar[d]^{f_{x}}\ar[r]&0\\
    W^{'}(S^{'}_{k})\ar[r]& W^{'}(x^{'}) \oplus W^{'}(S_{k-1})  \ar[r]&{W(x)}\ar[r]&0.
    }$$
      \end{enumerate}
 
    \end{enumerate}

 \begin{remark}
    By rephrasing the definition sequence of  map $V^{'}(x,y)$, we have the following commutative diagram
    
      \xymatrix{
    &&&&& V^{'}(x)\ar[rr]^{C_{V}(x)} \ar@{-->}[d]^{\exists !}\ar@(dl,ul)[dd]_{V^{'}(x, S_{k-1})}&& V(x^{'})\ar@{-->}[d]_{\exists !}\ar@(dr,ur)[dd]^{V(x^{'}, S_{k})}\\
    &&&&& V^{'}(y)\ar[d]^{V^{'}(y, S_{k-1})}\ar[rr]^{C_{V}(y)} &&V(y^{'}) \ar[d]_{V(y^{'}, S_{k})}\\
    &&&&& V(S_{k-1})\ar[rr]^{V(S_{k-1}, S_{k})}&&V(S_{k}).
    }

According to the uniqueness of map, there must exist a map $V^{'}(x,y): V^{'}(x) \to V^{'}(y)$, in a similar way, we have the map  $V(x^{'}, y^{'}): V(x^{'}) \to V(y^{'})$. Simplify the above diagram further, we have the final commutative diagram

     \xymatrix{
    &&&&& V^{'}(x)\ar[rr]^{C_{V}(x)}  \ar[d]^{V^{'}(x,y)}&& V(x^{'})\ar[d]^{V(x^{'}, y^{'})}\\
    &&&&& V^{'}(y)\ar[d]^{V^{'}(y, S_{k-1})}\ar[rr] ^{C_{V}(y)}&&V(y^{'}) \ar[d]^{V(y^{'}, S_{k})}\\
    &&&&& V(S_{k-1})\ar[rr]^{V(S_{k-1}, S_{k})}&&V(S_{k}).
    }
\end{remark}

\section {Main Results}
Let $A_\mathbb{R}= (\mathbb{R},\mathbf{S}, \preccurlyeq )$ be a continuous quiver of  type $A$, and $S_{k}$ be a sink. The new quiver  $\sigma _{k}A_{\mathbb{R}} = A^{'} _{\mathbb{R}} = (\mathbb{R}, {\mathbf{S}^{'}},\preccurlyeq^{'} )$ has been defined in Section 3. The category $rep(A_\mathbb{R})$ and the reflection functors $\mathsf{S}_{k}^{+}: rep(A_\mathbb{R})\longrightarrow rep(A^{'}_\mathbb{R})$, $\mathsf{S}_{k}^{-}: rep(A^{'}_\mathbb{R})\longrightarrow rep(A_\mathbb{R})$ have been defined in Section 3.

\begin{definition}
  Let $\overline{rep}(A_\mathbb{R})$ be the subcategory consisting of the following representations of $A_\mathbb{R}$,
  $$\overline{rep}(A_\mathbb{R})=\{ V \in rep(A_\mathbb{R}) | V(S_{k-1})\oplus V(S_{k+1})\overset{D_{1}^{'}}{\longrightarrow}V(S_{k})  , D_{1}^{'}\;  { is }\;{surjective}\} ,$$ where $D_{1}^{'}=(V(S_{k-1},S_{k}),- V(S_{k+1},S_{k}))$.
\end{definition}

\begin{definition}
Let $\underline{rep}(A^{'}_\mathbb{R})$ be the subcategory consisting of the following representations of ${rep}(A^{'}_\mathbb{R})$,
  $$\underline{rep}(A^{'}_\mathbb{R})=\{W\in rep(A_\mathbb{R}) | W(S^{'}_{k})\overset{D_{2}^{'}}{\longrightarrow} W(S_{k+1})\oplus W(S_{k-1}),D_{2}^{'} \; is\; injective \},$$ where $D_{2}^{'} =\binom{W(S^{'}_{k},S_{k+1})}{-W(S^{'}_{k},S_{k-1})} $.
\end{definition}

\begin{theorem}
The functor $\mathsf{S}_{k}^{+}: \overline{rep}(A_\mathbb{R}) \longrightarrow \underline{rep}(A^{'}_\mathbb{R})$ is an equivalence of categories.
\end{theorem}

For the proof of Theorem 4.3, we need the following lemmas.

\begin{lemma}
    With above notations, for any $x \in (S_{k-1}, S^{'}_{k} )$, the commutative diagram (5) is a pull-back
        \begin{equation}  \tag{5}
     \xymatrix{
    & V^{'}(S^{'}_{k})\ar[rr]^{V^{'}(S^{'}_{k}, S_{k+1})} \ar[d]^{V^{'}(S^{'}_{k}, x)} && V(S_{k+1})\ar[d]^{V(S_{k+1}, x^{'})} \\
    & V^{'}(x) \ar[rr]^{C_{V}(x)} &&{V(x^{'})},\\
    }
  \end{equation}
 and for any $x \in V(S^{'}_{k}, S_{k+1})$, the commutative diagram (6) is also a pull-back
  \begin{equation}  \tag{6}
     \xymatrix{
    & V^{'}(S^{'}_{k})\ar[rr]^{V^{'}(S^{'}_{k}, x)} \ar[d]^{V^{'}(x,S_{k-1} )\cdot V^{'}(S^{'}_{k}, x)} && V^{'}(x)\ar[d]^{C_{V}(x)} \\
    & V(S_{k-1}) \ar[rr] ^{^{V(S_{k-1}, x^{'})}}&&{V(x^{'})}.\\
    }
  \end{equation} 
   \begin{proof}
        Firstly, assume that $x \in (S_{k-1}, S^{'}_{k} )$, the commutative diagram (7) is a pull-back
          \begin{equation}\tag{7}
    \xymatrix{
    & V^{'}(x)\ar[d]^{V^{'}(x,S_{k-1} )} \ar[rr]^{C_{V}(x)} &&V(x^{'}) \ar[d]^{V(x^{'}, S_{k})}\\
    & V(S_{k-1})\ar[rr]^{V(S_{k-1} , S_{k})} &&{V(S_{k}),
    }}
        \end{equation}
  and the commutative diagram (8) is both a pull-back and a push-out       
          \begin{equation}  \tag{8}
   \xymatrix{
   & V^{'}(S^{'}_{k})\ar[rr]^{V^{'}(S^{'}_{k}, S_{k+1})} \ar[d]^{V^{'}(S^{'}_{k}, x)} && V(S_{k+1})\ar[d]^{V(S_{k+1}, x^{'})} \\
    & V^{'}(x)\ar[d]^{V^{'}(x,S_{k-1} )} \ar[rr]^{C_{V}(x)} &&V(x^{'}) \ar[d]^{V(x^{'}, S_{k})}\\
    & V(S_{k-1})\ar[rr]^{V(S_{k-1} , S_{k})} &&{V(S_{k}),
    }}
 \end{equation} 
 then we shall prove the commutative diagram (5) is a pull-back. That is to say, for $h: W \to V^{'}(x)$ and $k : W \to V(S_{k+1})$, if $C_{V}(x) \cdot h=V(S_{k+1}, x^{'}) \cdot k$, we need to find a homomorphism $l:W \rightarrow V^{'}(S^{'}_{k})$, such that $h=V^{'}(S^{'}_{k}, x) \cdot l$, $k=V^{'}(S^{'}_{k}, S_{k+1}) \cdot l$, as the following diagram shows
       
  \xymatrix{
    &&&& W\ar@/^/[rrrd]^{k} \ar@/_/[rdd]^{h} \ar@/_/[rddd]_{j} \ar@{.>}[rd]^l \\
    &&&&& V^{'}(S^{'}_{k})\ar[rr]^{V^{'}(S^{'}_{k}, S_{k+1})} \ar[d]^{V^{'}(S^{'}_{k}, x)} && V(S_{k+1})\ar[d]^{V(S_{k+1}, x^{'})} \\
    &&&&& V^{'}(x)\ar[d]^{V^{'}(x,S_{k-1} )} \ar[rr]^{C_{V}(x)} &&V(x^{'}) \ar[d]^{V(x^{'}, S_{k})}\\
    &&&&& V(S_{k-1})\ar[rr]^{V(S_{k-1} , S_{k})} &&V(S_{k}).
    }

Since the commutative diagram (8) is a pull-back, for homomorphism $V(S_{k-1}) \overset{m}{\longleftarrow} W \overset{k}{\longrightarrow} V(S_{k+1})$ such that $V(S_{k-1} , S_{k}) \cdot V^{'}(x,S_{k-1} ) \cdot V^{'}(S^{'}_{k}, x) \cdot l=V(x^{'}, S_{k}) \cdot V(S_{k+1}, x^{'}) \cdot k$, we can find a unique homomorphism $l:W \longrightarrow V^{'}(S^{'}_{k})$ satisfies the following diagram is commutative, that is,  $m=V^{'}(x,S_{k-1} ) \cdot V^{'}(S^{'}_{k}, x) \cdot l$, $k=V^{'}(S^{'}_{k}, S_{k+1}) \cdot l$, and $V(S_{k-1} , S_{k}) \cdot V^{'}(x,S_{k-1} ) \cdot V^{'}(S^{'}_{k}, x) \cdot l= V(x^{'}, S_{k}) \cdot V(S_{k+1}, x^{'}) \cdot k$, as the following diagram shows

    \xymatrix{
      &&&& W\ar@/^/[rrrd]^{k} \ar@/_/[rdd]_{V^{'}(x,S_{k-1} ) \cdot V^{'}(S^{'}_{k}, x) \cdot l} \ar[rd]^{l} \\
      &&&&& V^{'}(S^{'}_{k})\ar[rr]^{V^{'}(S^{'}_{k}, S_{k+1})} \ar[d]^{{V^{'}(x,S_{k-1} )} \cdot a} && V(S_{k+1})\ar[d]^{V(x^{'}, S_{k})\cdot V(S_{k+1}, x^{'})} \\
      &&&&& V(S_{k-1})\ar[rr]^{V(S_{k-1} , S_{k})} &&V(S_{k}).\\
      }
Since the commutative diagram (7) is a pull-back, for homomorphism $V(S_{k-1}) \overset{j}{\longleftarrow} W \overset{V(S_{k+1}, x^{'}) \cdot k}{\longrightarrow} V(x^{'})$ such that $V^{'}(x,S_{k-1} ) \cdot V^{'}(S^{'}_{k}, x) \cdot l=V(x^{'}, S_{k}) \cdot V(S_{k+1}, x^{'}) \cdot k$, we can find a unique homomorphism $h:W \longrightarrow V^{'}(x)$ satisfies the following diagram is commutative, that is, $j=V^{'}(x,S_{k-1} ) \cdot h$, $V(S_{k+1}, x^{'} ) \cdot k=C_{V}(x) \cdot h$, and $V^{'}(x,S_{k-1} ) \cdot h=V^{'}(x,S_{k-1} ) \cdot V^{'}(S^{'}_{k}, x) \cdot l$, as the following diagram shows
  
    \xymatrix{
      &&&& W\ar@/^/[rrrd]^{V(S_{k+1}, x^{'}) \cdot k} \ar@/_/[rdd]_{V^{'}(x,S_{k-1} ) \cdot V^{'}(S^{'}_{k}, x) \cdot l} \ar[rd]^{h} \\
      &&&&& V^{'}(x)\ar[rr]^{C_{V}(x)} \ar[d]^{V^{'}(x,S_{k-1} )} && V(x^{'})\ar[d]^{V(x^{'}, S_{k})}\\
      &&&&& V(S_{k-1})\ar[rr]^{V(S_{k-1} , S_{k})} &&V(S_{k}).
      }
    By the uniqueness of homomorphism $W \longrightarrow V(S_{k-1})$, we have $k=V^{'}(S^{'}_{k}, S_{k+1}) \cdot l$, $V^{'}(x,S_{k-1} ) \cdot V^{'}(S^{'}_{k}, x) \cdot l=V^{'}(x,S_{k-1} ) \cdot h$. Then we shall prove $V^{'}(S^{'}_{k}, x) \cdot l=h$.  
   
   For the homomorphism $h:W \rightarrow V^{'}(x)$, the following diagram is commutative
  
   \xymatrix{
        &&&& W\ar@/^/[rrrd]^{V(S_{k+1}, x^{'}) \cdot k} \ar@/_/[rdd]_{V^{'}(x,S_{k-1} ) \cdot h} \ar[rd]^{h} \\
        &&&&& V^{'}(x)\ar[rr]^{C_{V}(x)} \ar[d]^{V^{'}(x,S_{k-1} )} && V(x^{'})\ar[d]^{V(x^{'}, S_{k})} \\
        &&&&& V(S_{k-1})\ar[rr]^{V(S_{k-1} , S_{k})} &&V(S_{k}).
        }

  For the homomorphism $V^{'}(S^{'}_{k}, x) \cdot l:W \rightarrow V^{'}(x)$, we have $ C_{V}(x) \cdot V^{'}(S^{'}_{k}, x) \cdot l=V(S_{k+1}, x^{'}) \cdot V^{'}(S^{'}_{k}, S_{k+1}) \cdot l=V(S_{k+1}, x^{'}) \cdot k$, and the following diagram is commutative
  
   \xymatrix{
    &&&& W\ar@/^/[rrrd]^{V(S_{k+1}, x^{'}) \cdot k} \ar@/_/[rdd]_{V^{'}(x,S_{k-1} ) \cdot h} \ar[rd]^{V^{'}(S^{'}_{k}, x) \cdot l} \\
    &&&&& V^{'}(x)\ar[rr]^{C_{V}(x)} \ar[d]^{V^{'}(x,S_{k-1} )} && V(x^{'})\ar[d]^{V(x^{'}, S_{k})} \\
    &&&&& V(S_{k-1})\ar[rr]^{V(S_{k-1} , S_{k})} &&V(S_{k}).
    }

  By the uniqueness of map $W \rightarrow V^{'}(x)$, we have $h=V^{'}(S^{'}_{k}, x) \cdot l$, and $k=V^{'}(S^{'}_{k}, S_{k+1}) \cdot l$, that is, the commutative diagram (5) is a pull-back.  Similarly, the commutative diagram (6) is also a pull-back.
  
   \end{proof}
   
\end{lemma}

\begin{lemma}

With above notations, for any $x \in (S_{k-1}, S_{k} )$, the commutative diagram (9) is a push-out

  \begin{equation}  \tag{9}
    \xymatrix{
    & V^{'}(x^{'})\ar[d]^{V^{'}(x^{'} , S_{k-1})} \ar[rr]^{D_{V}(x)} &&\mathsf{S}_{k}^{-}V^{'}(x) \ar[d]^{\mathsf{S}_{k}^{-}V^{'}(x , S_{k})}\\
    & \mathsf{S}_{k}^{-}V^{'}(S_{k-1}) \ar[rr]^{\mathsf{S}_{k}^{-}V^{'}(S_{k-1}, S_{k})} &&{\mathsf{S}_{k}^{-}V^{'}(S_{k})
    },}
\end{equation} 
and for any $x \in (S_{k}, S_{k+1})$, the commutative diagram (10) is also a  push-out
 \begin{equation}  \tag{10}
    \xymatrix{
    & V^{'}(x^{'})\ar[d]^{D_{V}(x) }\ar[rr]^{V^{'}(x, S_{k+1}) } &&\mathsf{S}_{k}^{-}V^{'}(S_{k+1}) \ar[d]^{\mathsf{S}_{k}^{-}V^{'}(S_{k+1},  S_{k}) }\\
    & \mathsf{S}_{k}^{-}V^{'}(x)\ar[rr] ^{\mathsf{S}_{k}^{-}V^{'}(x , S_{k}) }&&{\mathsf{S}_{k}^{-}V^{'}(S_{k}).
    }}
\end{equation} 

   \begin{proof}
    Similarly to the proof of Lemma 4.4, we can get Lemma 4.5.
   \end{proof}
\end{lemma}

\begin{proposition}With above notations, we have $\mathsf{S}_{k}^{-}\mathsf{S}_{k}^{+} \cong 1_{\overline{rep}(A_\mathbb{R})}$, and $\mathsf{S}_{k}^{+}\mathsf{S}_{k}^{-} \cong 1_{\underline{rep}(A^{'}_\mathbb{R})}$.

 \begin{proof}

 (1) We shall first prove $\mathsf{S}_{k}^{-}\mathsf{S}_{k}^{+}V \cong V $.
 
   Assume that $x \in (S_{k},  S_{k+1})$. For any $V \in {rep}(A_\mathbb{R})$, $$(\mathsf{S}_{k}^{+}V(S^{'}_{k}), \mathsf{S}_{k}^{+}V(S^{'}_{k},S_{k+1} ), \mathsf{S}_{k}^{+}V(x^{'}, S_{k-1} ) \cdot \mathsf{S}_{k}^{+}V(S^{'}_{k}, x^{'} ))$$ is the pull-back of $$(V(S_{k}), V(S_{k+1}, x) \cdot V(x, S_{k}), V(S_{k-1}, S_{k})),$$ as the following diagram shows
 
    \xymatrix{
    &&& V(S_{k+1})\ar[d] ^{V(S_{k+1}, x)}  &&& \mathsf{S}_{k}^{+}V(S^{'}_{k})\ar[rr] ^{\mathsf{S}_{k}^{+}V(S^{'}_{k},S_{k+1} )}\ar[d]^{\mathsf{S}_{k}^{+}V(S^{'}_{k}, x^{'} )} && \mathsf{S}_{k}^{+}V(S_{k+1})\\
    &&&V(x) \ar[d]^{V(x, S_{k})} &  \ar[r]^{\mathsf{S}_{k}^{+}}&& \mathsf{S}_{k}^{+}V(x^{'}) \ar[d]^{\mathsf{S}_{k}^{+}V(x^{'}, S_{k-1} )}\\
    & V(S_{k-1})\ar[rr] ^{V(S_{k-1}, S_{k})}&&{V(S_{k})} &&& \mathsf{S}_{k}^{+}V(S_{k-1})
    .}

 For any $\mathsf{S}_{k}^{+}V \in {rep}^{'}(A_\mathbb{R})$, $(\mathsf{S}_{k}^{-}\mathsf{S}_{k}^{+}V(x), \mathsf{S}_{k}^{-}\mathsf{S}_{k}^{+}V(S_{k+1}, x) , D_{V}(x))$ is the push-out of $$(\mathsf{S}_{k}^{+}V(S^{'}_{k}) , \mathsf{S}_{k}^{+}V(S^{'}_{k},S_{k+1} ) , \mathsf{S}_{k}^{+}V(S^{'}_{k}, x^{'} )),$$
 and $({\mathsf{S}_{k}^{-} \mathsf{S}_{k}^{+}V(S_{k})} , \mathsf{S}_{k}^{+}V(x, S_{k} ) \cdot \mathsf{S}_{k}^{-}\mathsf{S}_{k}^{+}V(S_{k+1}, x) , \mathsf{S}_{k}^{-}\mathsf{S}_{k}^{+}V(S_{k+1}, S_{k}))$ is the push-out of $$(\mathsf{S}_{k}^{+}V(S^{'}_{k}) , \mathsf{S}_{k}^{+}V(S^{'}_{k},S_{k+1} ) , \mathsf{S}_{k}^{+}V(x^{'}, S_{k-1} ) \cdot \mathsf{S}_{k}^{+}V(S^{'}_{k}, x^{'} )),$$ as the following diagram shows 

   \xymatrix{
    \mathsf{S}_{k}^{+}V(S^{'}_{k})\ar[rr]^{\mathsf{S}_{k}^{+}V(S^{'}_{k},S_{k+1} ) } \ar[d]^{\mathsf{S}_{k}^{+}V(S^{'}_{k}, x^{'} )} && \mathsf{S}_{k}^{+}V(S_{k+1}) &&  \mathsf{S}_{k}^{+}V(S^{'}_{k})\ar@{-->}[rr] ^{\mathsf{S}_{k}^{+}V(S^{'}_{k},S_{k+1} ) } \ar@{-->}[d]^{\mathsf{S}_{k}^{+}V(S^{'}_{k}, x^{'} )} &&  \mathsf{S}_{k}^{-} \mathsf{S}_{k}^{+}V(S_{k+1})\ar[d]^{\mathsf{S}_{k}^{-}\mathsf{S}_{k}^{+}V(S_{k+1}, x)}\\
     \mathsf{S}_{k}^{+}V(x^{'}) \ar[d]^{\mathsf{S}_{k}^{+}V(x^{'}, S_{k-1} )} && \ar[r]^{\mathsf{S}_{k}^{-}}&& \mathsf{S}_{k}^{+}V(x^{'})\ar@{-->}[d]^{\mathsf{S}_{k}^{+}V(x^{'}, S_{k-1} )}  \ar@{-->}[rr] ^{D_{V}(x)} &&\mathsf{S}_{k}^{-} \mathsf{S}_{k}^{+}V(x) \ar[d]^{\mathsf{S}_{k}^{-}\mathsf{S}_{k}^{+}V(x, S_{k})}\\
     \mathsf{S}_{k}^{+}V(S_{k-1}) &&&&  \mathsf{S}_{k}^{-} \mathsf{S}_{k}^{+}V(S_{k-1})\ar[rr]^{\mathsf{S}_{k}^{-}\mathsf{S}_{k}^{+}V(S_{k+1}, S_{k})} &&{\mathsf{S}_{k}^{-} \mathsf{S}_{k}^{+}V(S_{k})}.} 
That is, by the definition of $\mathsf{S}_{k}^{-} \mathsf{S}_{k}^{+}V(x)$, the diagram (11) is a push-out

\begin{equation} \tag{11}
 \xymatrix{
     &  \mathsf{S}_{k}^{+}V(S^{'}_{k})\ar[rr] ^{\mathsf{S}_{k}^{+}V(S^{'}_{k},S_{k+1} )} \ar[d]^{\mathsf{S}_{k}^{+}V(S^{'}_{k}, x^{'} )}  &&  \mathsf{S}_{k}^{-} \mathsf{S}_{k}^{+}V(S_{k+1})\ar[d]^{\mathsf{S}_{k}^{-}\mathsf{S}_{k}^{+}V(S_{k+1}, x)}\\
     &\mathsf{S}_{k}^{+}V(x^{'})\  \ar[rr] ^{D_{V}(x)} &&{\mathsf{S}_{k}^{-} \mathsf{S}_{k}^{+}V(x) }
    .}
\end{equation}

By Lemma 4.4, the following diagram is a pull-back

  \xymatrix{
     &&&&  \mathsf{S}_{k}^{+}V(S^{'}_{k})\ar[rr] ^{\mathsf{S}_{k}^{+}V(S^{'}_{k},S_{k+1} )} \ar[d]^{\mathsf{S}_{k}^{+}V(S^{'}_{k}, x^{'} )} &&  V(S_{k+1})\ar[d]^{V(S_{k+1}, x)}\\
     &&&&\mathsf{S}_{k}^{+}V(x^{'})\  \ar[rr] ^{C_{V}(x^{'})}  &&{V(x) }
    .}
    
 Because $\mathsf{S}_{k}^{-} \mathsf{S}_{k}^{+}V(S_{k+1}) =  \mathsf{S}_{k}^{+}V(S_{k+1}) = V(S_{k+1})$, by the uniqueness of pull-back, $V(x) \cong \mathsf{S}_{k}^{-} \mathsf{S}_{k}^{+}V(x)$ for $x \in (S_{k}, S_{k+1})$. 
 
 If $x \in (S_{k-1} ,S_{k})$,  the conclusion is also true. That is, for any $x \in \mathbb{R}$, there exist a map $\eta _{V}(x): V(x) \to \mathsf{S}_{k}^{-} \mathsf{S}_{k}^{+}V(x)$ such that $V(x) \cong \mathsf{S}_{k}^{-} \mathsf{S}_{k}^{+}V(x)$.
 
 Next, we need to prove the diagram (12) is commutative
    \begin{equation}  \tag{12}
     \xymatrix{
    & V(x)\ar[rr] ^{\eta_{V}(x)}\ar[dd]^{V(x, y)}&& \mathsf{S}_{k}^{-}\mathsf{S}_{k}^{+}V(x)\ar[dd]^{(\mathsf{S}_{k}^{-}\mathsf{S}_{k}^{+}V)(x, y)}\\
    \\
    & V(y)\ar[rr] ^{\eta_{V}(y)}&&{\mathsf{S}_{k}^{-}\mathsf{S}_{k}^{+}V(y)
    }.} 
     \end{equation}
 
 For $y \prec x\in (S_{k},  S_{k+1})$, the following diagrams are commutative, respectively. By diagram (11), we have $\eta_{V}(y) \cdot V(S_{k+1}, y) = (\mathsf{S}_{k}^{-}\mathsf{S}_{k}^{+}V)(x , y) \cdot (\mathsf{S}_{k}^{-}\mathsf{S}_{k}^{+}V)(S_{k+1}, x)$. 
 
 \xymatrix{
    &&&& \mathsf{S}_{k}^{+}V(S^{'}_{k})\ar[d]_{\mathsf{S}_{k}^{+}V(S^{'}_{k}, x^{'} )} \ar[rr]^{\mathsf{S}_{k}^{+}V(S^{'}_{k},S_{k+1} )} && V(S_{k+1}) \ar[d]_{{V(S_{k+1}, x)}} \ar@(ur,dr)[rdd]^{V(S_{k+1}, y)} \\
    &&&& \mathsf{S}_{k}^{+}V(x^{'})\ar@/_/[rrrdd] _{D_{V}(x) \cdot (\mathsf{S}_{k}^{-}\mathsf{S}_{k}^{+}V)(x, y)}\ar[rr]^{C_{V}(x^{'})} &&V(x) \ar[rd]^{V(x , y)}\\
    &&&&&&& V(y) \ar[d]^{\eta_{V}(y)}\\
    &&&&&&& \mathsf{S}_{k}^{-}\mathsf{S}_{k}^{+}V(y) 
    }

 \xymatrix{
    &&&& \mathsf{S}_{k}^{+}V(S^{'}_{k})\ar[d]_{\mathsf{S}_{k}^{+}V(S^{'}_{k}, x^{'} )}  \ar[rr]^{\mathsf{S}_{k}^{+}V(S^{'}_{k},S_{k+1} )} && V(S_{k+1}) \ar[d]_{V(S_{k+1}, x)}\ar@(ur,dr)[rdd] ^{(\mathsf{S}_{k}^{-}\mathsf{S}_{k}^{+}V)(S_{k+1}, x)}\\
    &&&& \mathsf{S}_{k}^{+}V(x^{'})\ar@/_/[rrrdd]_{D_{V}(x) \cdot (\mathsf{S}_{k}^{-}\mathsf{S}_{k}^{+}V)(x, y)} \ar[rr]^{C_{V}(x^{'})} &&V(x)  \ar[rd]^{ \eta_{V}(x)}\\
    &&&&&&& {\mathsf{S}_{k}^{-}\mathsf{S}_{k}^{+}V(x) }\ar[d]^{(\mathsf{S}_{k}^{-}\mathsf{S}_{k}^{+}V)(x , y)}\\
    &&&&&&& \mathsf{S}_{k}^{-}\mathsf{S}_{k}^{+}V(y) 
    }

Since $( \mathsf{S}_{k}^{+}V(S^{'}_{k}) , \mathsf{S}_{k}^{+}V(S^{'}_{k},S_{k+1} ), \mathsf{S}_{k}^{+}V(S^{'}_{k}, x^{'} ))$ is the  push-out of $$(V(x), V(S_{k+1}, x), C_{V}(x^{'})),$$ we have $\eta_{V}(y) \cdot  V(x , y) = (\mathsf{S}_{k}^{-}\mathsf{S}_{k}^{+}V)(x , y) \cdot \eta_{V}(x)$, and the diagram (12) is commutative.

Similarly, for any $x \in (S_{k-1},  S_{k})$, we have the same conclusion.

 Thus we have $\mathsf{S}_{k}^{-}\mathsf{S}_{k}^{+}V \cong V $ for any $x \in (S_{k-1},  S_{k+1})$.

 (2) For maps $f_{x}: V(x) \to W(x) $ and $\mathsf{S}_{k}^{-}\mathsf{S}_{k}^{+}f_{x}: \mathsf{S}_{k}^{-}\mathsf{S}_{k}^{+}V(x) \to \mathsf{S}_{k}^{-}\mathsf{S}_{k}^{+}W(x)$, we need to prove the diagram (13) is commutative
   \begin{equation}  \tag{13}
    \xymatrix{
    & V(x)\ar[rr] ^{\eta_{V}(x)}\ar[dd]^{f_{x}}&& \mathsf{S}_{k}^{-}\mathsf{S}_{k}^{+}V(x)\ar[dd]^{\mathsf{S}_{k}^{-}\mathsf{S}_{k}^{+}f_{x}}\\
    \\
    & W(x)\ar[rr] ^{\eta_{W}(x)}&&{\mathsf{S}_{k}^{-}\mathsf{S}_{k}^{+}W(x)
    }.} 
   \end{equation}  
  
  Firstly, assume that $x \in (S_{k},  S_{k+1})$.
 For $V, W \in {rep}(A_\mathbb{R})$, the following diagrams are commutative, respectively. Note that $\eta_{W}(x) \cdot {f_x \cdot V(S_{k+1}, x)}=  \mathsf{S}_{k}^{-}\mathsf{S}_{k}^{+}f_x \cdot (\mathsf{S}_{k}^{-}\mathsf{S}_{k}^{+}V)(S_{k+1}, x) $.

 \xymatrix{
    &&&& \mathsf{S}_{k}^{+}V(S^{'}_{k})\ar[d]_{\mathsf{S}_{k}^{+}V(S^{'}_{k}, x^{'} )}  \ar[rr]^{\mathsf{S}_{k}^{+}V(S^{'}_{k},S_{k+1} )} && V(S_{k+1}) \ar[d]_{V(S_{k+1}, x)} \ar@(ur,dr)[rdd]^{f_x \cdot V(S_{k+1}, x)}\\
    &&&& \mathsf{S}_{k}^{+}V(x^{'})\ar@/_/[rrrdd]_ {D_{V}(x) \cdot \mathsf{S}_{k}^{-}\mathsf{S}_{k}^{+}f_{x} } \ar[rr]^{C_{V}(x^{'})} &&V(x) \ar[rd]^{f_x}\\
    &&&&&&& W(x) \ar[d]^{\eta_{W}(x)}\\
    &&&&&&& \mathsf{S}_{k}^{-}\mathsf{S}_{k}^{+}W(x)
    }

 \xymatrix{
    &&&& \mathsf{S}_{k}^{+}V(S^{'}_{k})\ar[d]_{\mathsf{S}_{k}^{+}V(S^{'}_{k}, x^{'} )}  \ar[rr]^{\mathsf{S}_{k}^{+}V(S^{'}_{k},S_{k+1} )} && V(S_{k+1}) \ar[d]_{V(S_{k+1}, x)} \ar@(ur,dr)[rdd] ^{(\mathsf{S}_{k}^{-}\mathsf{S}_{k}^{+}V)(S_{k+1}, x)}\\
    &&&& \mathsf{S}_{k}^{+}V(x^{'})\ar@/_/[rrrdd]_ {D_{V}(x) \cdot \mathsf{S}_{k}^{-}\mathsf{S}_{k}^{+}f_{x} } \ar[rr]^{C_{V}(x^{'})} &&V(x) \ar[rd]^{\eta_{V}(x)}\\
    &&&&&&& \mathsf{S}_{k}^{-}\mathsf{S}_{k}^{+}V(x) \ar[d]^{ \mathsf{S}_{k}^{-}\mathsf{S}_{k}^{+}f_x}\\
    &&&&&&& \mathsf{S}_{k}^{-}\mathsf{S}_{k}^{+}W(x)
    }

   Since $( \mathsf{S}_{k}^{+}V(S^{'}_{k}) , \mathsf{S}_{k}^{+}V(S^{'}_{k},S_{k+1} ) , \mathsf{S}_{k}^{+}V(S^{'}_{k}, x^{'} ) )$ is the  push-out of $$(V(x), V(S_{k+1}, x), C_{V}(x^{'})),$$we have $\eta_{W}(x) \cdot  f_x = \mathsf{S}_{k}^{-}\mathsf{S}_{k}^{+}f_x \cdot \eta_{V}(x)$, and the diagram (13) is commutative.
    
     If $x \in (S_{k-1}, S_{k})$, we have the same conclusion. By the arbitrariness of $x$, the following diagram is commutative
     
    \xymatrix{
    &&&&& V\ar[rr] ^{\eta_{V}}\ar[dd]_{f}&& \mathsf{S}_{k}^{-}\mathsf{S}_{k}^{+}V\ar[dd]^{\mathsf{S}_{k}^{-}\mathsf{S}_{k}^{+}f}\\
    \\
    &&&&& W\ar[rr] ^{\eta_{W}}&&{\mathsf{S}_{k}^{-}\mathsf{S}_{k}^{+}W
    }.}

 To sum up, we have $\mathsf{S}_{k}^{-}\mathsf{S}_{k}^{+} \cong 1_{\overline{rep}(A_\mathbb{R})}$. Similarly, we can get $\mathsf{S}_{k}^{+}\mathsf{S}_{k}^{-} \cong 1_{\underline{rep}(A^{'}_\mathbb{R})}$.
   \end{proof}
\end{proposition}

As the corollary of Proposition 4.6,  we get Theorem 4.3.

\bibliography{mybibfile}
\bibliographystyle{sci}


\end{document}